\documentclass[english]{smfart}
  \usepackage{macros_anglais}
  \usepackage{smfthm,smfenum}

\author{Beno\^{\i}t Kloeckner}
\address{Institut Fourier\\ 100 rue des Maths, BP 74\\ 38402 St Martin d'H\`eres, France}
\email{bkloeckn@fourier.ujf-grenoble.fr}
\urladdr{http://www-fourier.ujf-grenoble.fr/~bkloeckn/}

\title[Wasserstein space of Euclidean spaces]{A geometric study of Wasserstein spaces:\\ Euclidean spaces}

\theoremstyle{remark}
\newtheorem{ques}{Question}

\newcommand{\CAT}{\mathop {\mathrm{CAT}}\nolimits} 
\newcommand{\isom}{\mathop {\mathrm{Isom}}\nolimits}
\newcommand{\wass}{\mathop {\mathscr{W}}\nolimits}

\begin{document}

\begin{abstract}
In this article we consider Wasserstein spaces (with
quadratic transportation cost) as intrinsic
metric spaces. We are interested in usual geometric properties: curvature, rank
and isometry group, mostly in the case of Euclidean spaces.
Our most striking result is that the Wasserstein space
of the line admits ``exotic'' isometries, which
do not preserve the shape of measures.
\end{abstract}


\maketitle

\setcounter{tocdepth}{1}
\setcounter{secnumdepth}{2}
\tableofcontents

\section{Introduction}

The concept of optimal transportation recently raised a growing interest
in links with the geometry of metric spaces. In particular the $L^2$ 
Wasserstein space $\wass_2(X)$ have been used by Von Renesse and Sturm \cite{vonRenesse-Sturm},
Sturm \cite{Sturm} and Lott and Villani \cite{Villani-Lott} to define certain curvature conditions
on a metric space $X$.
Many useful properties are inherited from $X$ by $\wass_2(X)$ (separability, completeness, geodesicness, 
some non-negative curvature conditions) while some other
are not, like local compacity.

In this paper, we study the geometry of Wasserstein spaces
as intrinsic spaces. We are interested, for example, in the isometry group
of $\wass_2(X)$, in its curvature and in its rank (the greatest possible dimension of 
a Euclidean space that embeds in it). In the case of the Wasserstein space
of a Riemannian manifold, itself seen as an infinite-dimensional Riemannian manifold,
the Riemannian connection and curvature have been computed by Lott \cite{Lott}.
See also \cite{Takatsu} where Takatsu studies the subspace of Gaussian measures 
in $\wass_2(\mR^n)$, and \cite{Ambrosio-Gigli} where Ambrosio and Gigli are interested
in the second order analysis on $\wass_2(\mR^n)$, in particular its parallel transport.

The Wasserstein space $\wass_2(X)$ contains a copy of $X$, the image of the isometric
embedding 
\begin{eqnarray*}
E : X &\to& \wass_2(X) \\
x &\mapsto& \delta_x
\end{eqnarray*}
where $\delta_x$ is the Dirac mass at $x$.
Moreover, given an isometry $\varphi$ of $X$ one defines an isometry
$\varphi_\#$ of $\wass_2(X)$ by $\varphi_\#(\mu)(A)=\mu(\varphi^{-1}(A))$.
We thus get an embedding
$$\#:\isom X\to\isom \wass_2(X)$$
These two elementary facts connect the geometry of $\wass_2(X)$ to that of $X$.

One could expect that $\#$
is onto, \textit{i.e.} that all isometries of $\wass_2(X)$ are induced
by those of $X$ itself. Elements of $\#(\isom X)$ are called \emph{trivial}
isometries. Let us introduce a weaker property: a self-map
$\Phi$ of $\wass_2(X)$ is said to \emph{preserve shapes} if for all
$\mu\in\wass_2(X)$, there is an isometry $\varphi$ of $X$ (that depends upon $\mu$) such that
$\Phi(\mu)=\varphi_\#(\mu)$. An isometry that does not
preserve shapes is said to be \emph{exotic}.

Our main result is the surprising fact that $\wass_2(\mR)$ admits
exotic isometries. More precisely we prove the following.
\begin{theo}\label{enonce:isometries}
The isometry group of $\wass_2(\mR)$ is a semidirect product
\begin{equation}
\isom\mR\ltimes\isom\mR
\label{eq:isometries}
\end{equation}
Both factors decompose: $\isom\mR=\mZ/2\mZ\ltimes\mR$
and the action defining the semidirect product \eqref{eq:isometries} is simply given by the usual
action of the left $\mZ/2\mZ$ factor on the right $\mR$ factor, that is
$(\varepsilon,v)\cdot(\eta,t)=(\eta,\varepsilon t)$ where $\mZ/2\mZ$ is
identified with $\{\pm1\}$.

In \eqref{eq:isometries}, the left factor is the image of $\#$ and
the right factor consist in all isometries that fix pointwise the
set of Dirac masses. In the decomposition of the latter, 
the $\mZ/2\mZ$ factor is generated by a non-trivial
involution that preserves shapes, while the $\mR$ factor
is a flow of exotic isometries.
\end{theo}

The main tool we use is the explicit description of the geodesic between two points
$\mu_0$, $\mu_1$ of $\wass_2(\mR)$ that follows from the fact that the unique optimal transportation
plan between $\mu_0$ and $\mu_1$ is the non-decreasing rearrangement. It implies
that most of the geodesics in $\wass_2(\mR)$ are not complete, and we rely on this fact
to give a metric characterization of Dirac masses
and of linear combinations of two Dirac masses, among all points of $\wass_2(X)$. We also
use the fact that $\wass_2(\mR)$ has vanishing curvature in the sense of Alexandrov.

Let us describe roughly the non-trivial isometries that fix pointwise the set of
Dirac masses.
On the one hand, the non-trivial isometry generating the $\mZ/2\mZ$ factor is defined as follows: a measure $\mu$ is
mapped to its symmetric with respect to its center of mass.
On the other hand, the exotic isometric flow tends to put all the mass on one side of the center
of gravity (that must be preserved), close to it, and to send a small bit of mass far away on the other
side (so that the Wasserstein distance to the center of mass is preserved). In particular,
under this flow any measure $\mu$ converges weakly (but of course not in $\wass_2(\mR)$)
to $\delta_x$ (where $x$ is the center of mass
of $\mu$), see Proposition \ref{enonce:weak_convergence}. 

The case of the line seems very special. For example,
$\wass_2(\mR^n)$ admits non-trivial isometries but
all of them preserve shapes.
\begin{theo}\label{enonce:isom_highdim}
If $n\geqslant 2$, the isometry group of $\wass_2(\mR^n)$ is a semidirect product
\begin{equation}
\isom(\mR^n)\ltimes\mathrm{O}(n)
\end{equation}
where the action of an element $\psi\in\isom(\mR^n)$ on $\mathrm{O}(n)$
is the conjugacy by its linear part $\vec{\psi}$.

The left factor is the image of $\#$
and each element in the right factor fixes all Dirac masses and preserves shapes.
\end{theo}
The proof relies on Theorem \ref{enonce:isometries}, some elementary properties
of $L^2$ optimal transportation in $\mR^n$ and Radon's Theorem \cite{Radon}.

We see that the quotient $\isom \wass_2(\mR^n)/\isom\mR^n$ is compact if and only
if $n>1$. The higher-dimensional Euclidean spaces
are more rigid than the line for this problem, and we expect  most of the other metric spaces
to be even more rigid in the sense that $\#$ is onto.

Another consequence of the study of complete geodesics concerns the rank of
$\wass_2(\mR^n)$.
\begin{theo}\label{enonce:flats}
There is no isometric embedding of $\mR^{n+1}$ into $W_2(\mR^n)$.
\end{theo}
 It is simple to prove that despite Theorem \ref{enonce:flats}, large pieces of
$\mR^n$ can be embedded into $\wass_2(\mR)$, which has consequently
infinite weak rank in a sense to be precised.
 As a consequence, we get for example:
\begin{prop}\label{enonce:not_hyperbolic}
If $X$ is any Polish geodesic metric space that contains a complete geodesic,
then $\wass_2(X)$ is not $\delta$-hyperbolic.
\end{prop}
This is not surprising, since it is well-known that the negative curvature assumptions tend
not to be inherited from $X$ by its Wasserstein space. An explicit example is computed in 
\cite{Ambrosio-Gigli-Savare} (Example 7.3.3); more generaly, if $X$ contains a rhombus
(four distinct points $x_1,x_2,x_3,x_4$ so that $d(x_i,x_{i+1})$ is independent of the cyclic index $i$)
then $\wass_2(X)$ is not uniquely geodesic, and in particular not $\CAT(0)$, even if $X$ itself is
strongly negatively curved. 

\subsection*{Organization of the paper}

Sections \ref{sec:recalls} to \ref{sec:curvature} collect some properties
needed in the sequel.
Theorem \ref{enonce:isometries} is proved in Section \ref{sec:isometries},
Theorem \ref{enonce:isom_highdim} in Section \ref{sec:plane}.
Section \ref{sec:ranks} is devoted to the ranks of $\wass_2(\mR)$ and
$\wass_2(\mR^n)$, and we end in Section \ref{sec:open} with some open questions.

\subsection*{Acknowledgements}

I wish to thank all speakers of the workshop on optimal transportation held in the
Institut Fourier in Grenoble, especially Nicolas Juillet with whom I had numerous discussion
on Wasserstein spaces, and its organizer Herv\'e Pajot. I am also indebted to Yann
Ollivier for advises and pointing out some inaccuracies and mistakes in preliminary
versions of this paper.

\section{The Wasserstein space}\label{sec:recalls}

In this preliminary section we recall well-known general facts on $\wass_2(X)$. One can refer
to \cite{Villani,Villani2} for further details and much more.
Note that the denomination ``Wasserstein space'' is debated and
historically inaccurate. However, it is now the most common denomination
and thus an occurrence of the self-applying theorem of Arnol'd according to which a mathematical
result or object is usually attributed to someone that had little to do with it.

\subsection{Geodesic spaces}

Let $X$ be a Polish (\textit{i.e.} complete and separable metric) space, and assume that $X$
is geodesic, that is: between two points there is a rectifiable curve
whose length is the distance between the considered points. Note that we only
consider \emph{globally} minimizing geodesics, and that a geodesic is always
assumed to be parametrized proportionally to arc length.

One defines the Wasserstein space of $X$ as the set $\wass_2(X)$ of Borel probability measures
$\mu$ on $X$ that satisfy
\[\int_X d^2(x_0,x) \mu(dx)<+\infty\]
for some (hence all) point $x_0\in X$, equipped by the distance $d_{\wass}$ defined by:
\[d_{\wass}^2(\mu_0,\mu_1)=\inf\int_{X\times X} d^2(x,y) \Pi(dxdy)\]
where the infimum is taken over all coupling $\Pi$ of $\mu_0$, $\mu_1$. A coupling
realizing this infimum is said to be \emph{optimal}, and there always exists an optimal coupling.

The idea behind this distance is linked to the Monge-Kantorovitch problem: given a
unit quantity of goods distributed in $X$ according to $\mu_0$, what is the most
economical way to displace them so that they end up distributed according to $\mu_1$,
when the cost to move a unit of good from $x$ to $y$ is given by $d^2(x,y)$? The minimal cost
is $d_{\wass}^2(\mu_0,\mu_1)$ and a transportation plan achieving this minimum is
an optimal coupling.

An optimal coupling is said to be \emph{deterministic} if it can be written under the form
$\Pi(dx dy)=\mu(dx)\mathds{1}[y=Tx]$ where $T:X\to X$ is a measurable map and $\mathds{1}[A]$ is $1$ 
if $A$ is satisfied and $0$ otherwise. This means that the coupling does not split
mass: all the mass at point $x$ is moved to the point $Tx$. One usually write $\Pi=(\id\times T)_\#\mu$.
Of course, for $\Pi$ to be a coupling between $\mu$ and $\nu$, the relation
$\nu=T_\#\mu$ must hold.

Under the assumptions we put on $X$, the metric space $\wass_2(X)$ is itself Polish and geodesic.
If moreover $X$ is uniquely geodesic, then to each optimal coupling $\Pi$ between $\mu_0$ and $\mu_1$
is associated a unique geodesic in $\wass_2(X)$ in the following way. Let $C([0,1],X)$ be
the set of continuous curves $[0,1]\to X$, let $g:X\times X\to C([0,1],X)$
be the application that maps $(x,y)$ to the constant speed geodesic between these points, and
for each $t\in[0,1]$ let $e(t):C([0,1],X)\to X$ be the map $\gamma\mapsto \gamma(t)$.
Then $t\mapsto e(t)_\# g_\# \Pi$ is a geodesic between $\mu_0$ and $\mu_1$. Informally, this means
that we choose randomly a couple $(x,y)$ according to the joint law $\Pi$,
then take the time $t$ of the geodesic $g(x,y)$. This gives a random point in $X$, whose
law is $\mu_t$, the time $t$ of the geodesic in $\wass_2(X)$ associated to the optimal coupling $\Pi$.
Moreover, all geodesics are obtained that way.

Note that for most spaces $X$, the optimal coupling is not unique for all pairs of probability measures, and
$\wass_2(X)$ is therefore not uniquely geodesic even if $X$ is.

One of our goal
is to determine whether the Dirac measures can be detected inside $\wass_2(X)$ by purely geometric properties,
so that we can link the isometries of $\wass_2(X)$ to those of $X$.

\subsection{The line}

Given the distribution function 
\[F:x\mapsto \mu(]-\infty,x])\]
of a probability measure $\mu$, one defines its left-continuous inverse:
\begin{eqnarray*}
F^{-1} : \left] 0,1\right[ &\to& \mR \\
         m   &\mapsto&\sup\{x\in\mR ; F(x)\leqslant m\}
\end{eqnarray*}
that is a non-decreasing, left-continuous function; $\lim_0 F^{-1}$ is the infimum of the support
of $\mu$ and $\lim_1 F^{-1}$ its supremum. A discontinuity of $F^{-1}$ happens for each interval
that does not intersect
the support of $\mu$, and $F^{-1}$ is constant on an interval for each atom of $\mu$.

\begin{figure}[ht]\begin{center}
\input{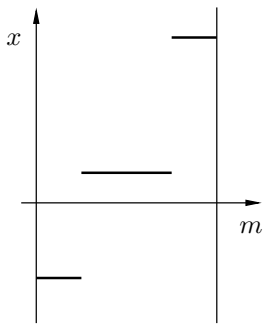}
\caption{Inverse distribution function of a combination of three Dirac masses}
\end{center}\end{figure}

Let $\mu_0$ and $\mu_1$ be two points of $\wass_2(\mR)$, and let $F_0,F_1$ be their distribution functions.
Then the distance between $\mu_0$ and $\mu_1$ is given by
\begin{equation}
d^2(\mu_0,\mu_1)=\int_0^1 \left(F_0^{-1}(m)-F_1^{-1}(m)\right)^2dm \label{eq:distance}
\end{equation}
and there is a unique constant speed geodesic $(\mu_t)_{t\in[0,1]}$, where $\mu_t$ has a distribution
function $F_t$ defined by
\begin{equation}
F_t^{-1}=(1-t)F_0^{-1}+t F_1^{-1}\label{eq:geodesic_desc}
\end{equation}
This means that the best way to go from $\mu_0$ to $\mu_1$ is simply to rearrange increasingly the mass,
a consequence of the convexity of the cost function. For example, if $\mu_0$ and $\mu_1$ are uniform 
measures on $[0,1]$ and $[\varepsilon,1+\varepsilon]$, then the optimal coupling is deterministic
given by the translation $x\mapsto x+\varepsilon$. That is: the best way to go from $\mu_0$ to $\mu_1$
is to shift every bit of mass by $\varepsilon$. If the cost function where linear, it would be equivalent
to leave the mass on $[\varepsilon, 1]$ where it is and move the remainder from $[0,\varepsilon]$ to
$[1,1+\varepsilon]$. If the cost function where concave, then the latter solution would be better than the former.

\begin{figure}[ht]\begin{center}
\input{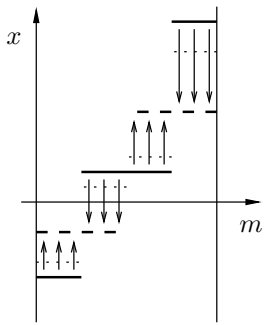}
\caption{A geodesic between two atomic measures: the mass moves with speed proportional to the
         length of the arrows.}
\end{center}\end{figure}

\subsection{Higher dimensional Euclidean spaces}\label{sec:Rn_bases}

The Monge-Kantorovich problem is far more intricate in 
$\mR^n$ ($n\geqslant 2$) than in $\mR$. The major contributions
of Knott and Smith \cite{KS1,KS2} and Brenier \cite{Brenier1,Brenier2}
give a quite satisfactory characterization of optimal couplings
and their unicity when the two considered measures
$\mu$ and $\nu$ are absolutely continuous (with respect
to the Lebesgue measure). We shall not give details of these works,
for which we refer to \cite{Villani, Villani2} again. Let us however
consider some toy cases, which will prove useful later on.
Missing proofs can be found in \cite{Juillet}, Section 2.1.2.

We consider $\mR^n$ endowed with its canonical inner product
and norm, denoted by $|\cdot|$.

\subsubsection{Translations}

Let $T_v$ be the translation of vector $v$ and assume
that $\nu=(T_v)_\# \mu$. Then the unique optimal coupling
between $\mu$ and $\nu$ is deterministic, equal
to $(\id\times T_v)_\#\mu$, and therefore $d_{\wass}(\mu,\nu)=|v|$.
This means that the only
most economic way to move the mass from $\mu$ to $\nu$
is to translate each bit of mass by the vector $v$. This is a quite
intuitive consequence of the convexity of the cost.
In particular, the geodesic between $\mu$ and $\nu$ can be extended
for all times $t\in\mR$. This happens only in this case as we shall
see later on.

\subsubsection{Dilations}

Let $D_x^\lambda$ be the dilation of center $x$ and ratio $\lambda$
and assume that $\nu=(D_x^\lambda)_\# \mu$. Then the unique
optimal coupling between $\mu$ and $\nu$ is deterministic,
equal to $(\id\times D_x^\lambda)_\# \mu$. In particular,
$$d_{\wass}(\mu,\nu)=|1-\lambda^2|^{\frac12}d_{\wass}(\mu,\delta_x)$$

As a consequence, the geodesic between $\mu$ and $\nu$ is unique and made
of homothetic of $\mu$, and can be extended only to a semi-infinite interval: 
it cannot be extended beyond $\delta_x$ (unless $\mu$ is a Dirac mass itself).

\subsubsection{Orthogonal measures}

Assume that $\mu$ and $\nu$ are supported on orthogonal affine subspaces 
$V$ and $W$ of $\mR^n$. 
Then if $\Pi$ is any coupling, assuming $0\in V\cap W$, we have
\begin{eqnarray*}
\int_{\mR^n\times\mR^n} |x-y|^2 \Pi(dx dy)
  &=&\int_{\mR^n\times\mR^n} (|x|^2 + |y|^2)\Pi(dx dy)\\
  &=&\int_V |x|^2\mu(dx)+\int_W |y|^2\nu(dy)
\end{eqnarray*}
therefore the cost is the same whatever the coupling.

\subsubsection{Balanced combinations of two Dirac masses}

\begin{figure}[htp]\begin{center}
\includegraphics{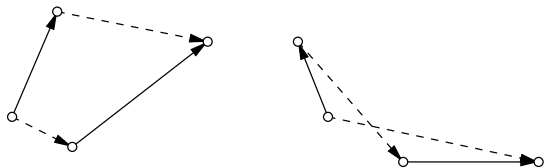}
\caption{Optimal coupling between balanced combinations of two Dirac masses. 
         Continuous arrows represent the vectors $y_0-x_0$ and $y_1-x_1$ while dashed arrows
         represent the optimal coupling.}\label{fig:transport_ex1}
\end{center}\end{figure}

Assume that $\mu=1/2\delta_{x_0}+1/2\delta_{y_0}$ and
$\nu=1/2\delta_{x_1}+1/2\delta_{y_1}$. A coupling between
$\mu$ and $\nu$ is entirely determined by the amount $m\in[0,1/2]$ of
mass sent from  $x_0$ to $x_1$. The cost of the coupling is
$$\frac12 |x_1-y_0|^2+\frac12|x_0-y_1|^2-2m(y_0-x_0)\cdot(y_1-x_1)$$
thus the optimal coupling is unique and deterministic if $(y_0-x_0)\cdot(y_1-x_1)\neq0$,
given by the map $(x_0,y_0)\mapsto (x_1,y_1)$ if $(y_0-x_0)\cdot(y_1-x_1)>0$
and by the map $(x_0,y_0)\mapsto (y_1,x_1)$ if $(y_0-x_0)\cdot(y_1-x_1)<0$ (figure \ref{fig:transport_ex1}).
Of course if $(y_0-x_0)\cdot(y_1-x_1)=0$, then all coupling have the same cost and
are therefore optimal.

If the combinations are not balanced (the mass is not equally split between the two point
of the support), then the optimal coupling is easy to deduce from the preceding computation.
For example if $(y_0-x_0)\cdot(y_1-x_1)>0$ then as much mass as possible must be sent
from $x_0$ to $x_1$, and this determines the optimal coupling.

This example has a much more general impact than it might seem: it can be generalized to
the following (very) special case of the \emph{cyclical monotonicity} (see for example
\cite{Villani2}, Chapter 5) which will prove useful in the sequel.
\begin{lemm}\label{lemm:cyclical_monotonicity}
If $\Pi$ is an optimal coupling between any two probability measures on $\mR^n$,
then
$$(y_0-x_0)\cdot(y_1-x_1)\geqslant 0$$
holds whenever $(x_0,x_1)$ and $(y_0,y_1)$ are in the support of $\Pi$.
\end{lemm}

\subsection{Spaces of nonpositive curvature}

In this paper we shall consider two curvature conditions. The first one is a negative curvature 
condition, the $\delta$-hyperbolicity introduced by Gromov (see for example \cite{Bridson-Haeffliger}).
 A geodesic space is said to be $\delta$-hyperbolic (where $\delta$
is a non-negative number) if in any triangle, any point of any of the sides is at distance
at most $\delta$ from one of the other two sides. For example, the real hyperbolic space is
$\delta$-hyperbolic (the value of $\delta$ depending on the value of the curvature),
a tree is $0$-hyperbolic and the euclidean spaces of dimension at least $2$
are not $\delta$-hyperbolic for any $\delta$.

The second condition is the classical non-positive sectional curvature condition $\CAT(0)$, detailed
in Section \ref{sec:curvature}, that roughly means that triangles are thinner in $X$ than in 
the euclidean plane.
Euclidean spaces, any Riemannian manifold having non-positive sectional curvature are examples
of locally $\CAT(0)$ spaces.

A geodesic $\CAT(0)$ Polish space $X$ is also called a \emph{Hadamard space}. 
A Hadamard space is uniquely geodesic, and admits a natural boundary at infinity.
The feature that interests us most is the following classical result: if $X$ is a Hadamard space,
given $\mu\in \wass_2(X)$ there is a unique point $x_0\in X$, called the center of mass
of $\mu$, that minimizes the quantity $\int_X d^2(x_0,x)\mu(dx)$.
If $X$ is $\mR^n$ endowed with the canonical scalar product, then the center of mass is of course
$\int_{\mR^n} x \mu(dx)$ but in the general case, the lack of an affine structure on $X$
prevents the use of such a formula.

We thus get a map $P: \wass_2(X)\to X$ that maps any $L^2$ probability measure to its center of mass.
Obviously, $P$ is a left inverse to $E$ and one can hope to use this map to link closer
the geometry of $\wass_2(X)$ to that of $X$. That's why our questions, unlike most of the classical ones
in optimal transportation, might behave more nicely when the curvature is non-positive than when it is non-negative.

\section{Geodesics}\label{sec:geodesics}

The content of this section, although difficult to locate in the bibliography, is part of the folklore
and does not pretend to originality. We give proofs for the sake of completeness.

\subsection{Case of the line}

We now consider the geo\-de\-sics of $\wass_2(\mR)$. Our first
goal is to determine on which maximal interval they can be
extended.

\subsubsection{Maximal extension}

Let $\mu_0$, $\mu_1$ be two points of $\wass_2(\mR)$ and $F_0$, $F_1$ their distribution functions.
Let $(\mu_t)_{t\in[0,1]}$ be the geodesic between $\mu_0$ and $\mu_1$. Since $\wass_2(\mR)$
is uniquely geodesic, there is a unique maximal interval on which $\gamma$ can be extended into
a geodesic, denoted by $I(\mu_0,\mu_1)$. 

\begin{lemm}\label{enonce:max_extension}
One has
\[I(\mu_0,\mu_1)=\{t\in\mR ; F_t^{-1} \mbox{ is non-decreasing}\}\]
where $F_t^{-1}$ is defined by the formula \eqref{eq:geodesic_desc}.
It is a closed interval. If one of its bound $t_0$ is finite,
then $\mu_{t_0}$ does not have bounded density with respect to the Lebesgue
measure.
\end{lemm}

\begin{proof}
Any non-decreasing left continuous function is the inverse distribution
function of some probability measure. If such a function is obtained
by an affine combination of probabilities belonging to $\wass_2(\mR)$,
then its probability measure belongs to $\wass_2(\mR)$ too.

Moreover, an affine combination of two left continuous function
is left continuous, so that 
\[I(\mu_0,\mu_1)=\{t\in\mR ; F_t^{-1} \mbox{ is non-decreasing}\}.\]

The fact that $I(\mu_0,\mu_1)$ is closed follows from the stability
of non-decreasing functions under pointwise convergence.

If the minimal slope
\[\inf\left\{\frac{F_t^{-1}(m)-F_t^{-1}(m')}{m-m'} ; 0<m<m'<1\right\}\]
is positive for some $t$, then it stays positive in a neighborhood of $t$.
Thus, a finite bound of $I(\mu_0,\mu_1)$ must have zero minimal slope,
and cannot have a bounded density.
\end{proof}

A geodesic is said to be complete if it is defined for all times.
We also consider geodesic rays, defined on an interval $[0,T]$ or $[0,+\infty[$
(in the latter case we say that the ray is complete),
and geodesic segments, defined on a closed interval.

It is easy to deduce a number of consequences from Lemma \ref{enonce:max_extension}.

\begin{prop}\label{enonce:consequence} In $\wass_2(\mR)$:
\begin{enumerate}
\item\label{it:consequence_a} any geodesic ray issued from a Dirac mass can be extended to a complete
      ray,
\item\label{it:consequence_b} no geodesic ray issued from a Dirac mass can be extended for negative times,
      except if all of its points are Dirac masses,
\item\label{it:consequence_c} up to normalizing the speed, the only complete geodesics 
      are those obtained by translating a point of $\wass_2(\mR)$:
      \[\mu_t(A)=\mu_0(A-t),\]
\end{enumerate}
\end{prop}

\begin{proof}
The inverse distribution function of a Dirac mass $\delta_x$ is the constant function
$F_0^{-1}$ with value $x$. Since it slopes 
\[\frac{F_0^{-1}(m)-F_0^{-1}(m')}{m-m'}\qquad0<m<m'<1\]
are all zero, for all positive times $t$ the functions $F_t^{-1}$ defined by 
formula \eqref{eq:geodesic_desc}
for any non-decreasing $F_1^{-1}$ are non-decreasing. However, for $t<0$ the $F_t^{-1}$ 
are not non-decreasing if $F_1^{-1}$ is not constant, 
we thus get (\ref{it:consequence_a}) and (\ref{it:consequence_b}).

Consider a point $\mu_0$ of $\wass_2(\mR)$ defined by an inverse distribution function
$F_0^{-1}$, and consider a complete geodesic $(\mu_t)$ issued from $\mu_0$.
Let $F_t^{-1}$ be the inverse distribution function of $\mu_t$. Then, since
$\mu_t$ is defined for all times $t>0$, the slopes of $F_1^{-1}$ must be greater
than those of $F_0^{-1}$:
\[F_0^{-1}(m)-F_0^{-1}(m')\leqslant F_1^{-1}(m)-F_1^{-1}(m')\quad\forall m < m'\]
otherwise, when $t$ increases, some slope of $F_t^{-1}$ will decrease linearly
in $t$, thus becoming negative in finite time.

But since $\mu_t$ is also defined for all $t<0$, the slopes of $F_1^{-1}$ must be
lesser than those of $F_0^{-1}$. They are therefore equal, and the two inverse
distribution function are equal up to an additive constant. The geodesic $\mu_t$
is the translation of $\mu_0$ and we proved (\ref{it:consequence_c}).
\end{proof}

\subsubsection{Convex hulls of totally atomic measures}

Define in $\wass_2(\mR)$ the following sets:
\begin{eqnarray*}
\Delta_1&=&\{\delta_x ; x\in \mR\} \\
\Delta_n&=&\{\sum_{i=1}^n a_i\delta_{x_i} ; x_i\in\mR, \textstyle\sum a_i = 1\} \\
\Delta'_{n+1} &=&\Delta_{n+1}\setminus\Delta_n
\end{eqnarray*}

Recall that if $X$ is a Polish geodesic space and $C$ is a subset of $X$,
one says that $C$ is convex if every geodesic segment whose endpoints
are in $C$ lies entirely in $C$.
The convex hull of a subset $Y$ is the least convex set $C(Y)$ that contains $Y$. 
It is well defined since the intersection of two convex sets is a convex set,
and is equal to $\cup_{n\in\mN} Y_n$ where $Y_0=Y$ and $Y_{n+1}$ is obtained
by adding to $Y_n$ all points lying on a geodesic whose endpoints are in $Y_n$.

Since $\Delta_1$ is the image of the isometric embedding $E:\mR\to \wass_2(\mR)$,
it is a convex set. This is not the case of $\Delta_n$ is $n>1$.
In fact, we have the following.

\begin{prop}\label{enonce:convex}
If $n>1$, any point $\mu$ of $\Delta_{n+1}$ lies on a geodesic segment
with endpoints in $\Delta_n$. Moreover, the endpoints can be chosen with
the same center of mass than that of $\mu$. 
\end{prop}

\begin{proof}
If $\mu\in\Delta_n$ the result is obvious.
Assume $\mu=\sum a_i \delta_{x_i}$ is in $\Delta'_{n+1}$.
We can assume further that $x_1<x_2<\cdots< x_{n+1}$. Consider the measures
\begin{eqnarray*}
\mu_{-1}   &=& \sum_{i<n-1} a_i\delta_{x_i}+(a_{n-1}+a_n)\delta_{x_{n-1}}
            +a_{n+1}\delta_{x_{n+1}} \\
\mu_1   &=& \sum_{i<n-1} a_i\delta_{x_i}+a_{n-1}\delta_{x_{n-1}}+(a_n+a_{n+1})\delta_{x_{n+1}}.
\end{eqnarray*}
Then $\mu$ lies on the geodesic segment from $\mu_{-1}$ to $\mu_1$.
To get a constant center of mass, one considers the geodesic
\[\mu_t = \sum_{i\le n-1} a_i\delta_{x_i}+a_{n}\delta_{x_{n}+t}+a_{n+1}\delta_{x_{n+1}-\alpha t}\]
where $\alpha=\frac{a_{n}}{a_{n+1}}$.
\end{proof}

In particular, we get the
following noteworthy fact that will
prove useful latter on.

\begin{prop}\label{enonce:convex_hull}
The convex hull of $\Delta_n$ is dense in $\wass_2(\mR)$ if $n>1$.
\end{prop}

\begin{proof}
Follows from Proposition \ref{enonce:convex} since the set of totally ato\-mic measures
$\bigcup_n \Delta_n$ is dense in $\wass_2(\mR)$.
\end{proof}

\subsection{Complete geodesics in higher dimension}\label{sec:geod_high}

In $\mR^n$, the optimal coupling and thus the geodesics
are not as explicit as in the case of the line. 
It is however
possible to determine which geodesic can be extended to all
times in $\mR$.

\begin{lemm}
Let $\mu=(\mu_t)_{t\in I}$ be a geodesic in $\wass_2(\mR^n)$
associated to an optimal coupling $\Pi$ between $\mu_0$ and
$\mu_1$. Then for all times $r$ and $s$ in $I$ and all
pair of points $(x_0,x_1),(y_0,y_1)$ in the support
of $\Pi$, the following hold:
$$|u|^2+(r+s)u\cdot v+rs|v|^2\geqslant0$$
where $u=y_0-x_0$ and $v=y_1-x_1-(y_0-x_0)$.
\end{lemm}

\begin{proof}
Let us introduce the following notations: for all pair of points
$a_0,a_1\in\mR^n$, $a_t=(1-t)a_0+ta_1$
and $\Pi_{r,s}$ is the law of the random variable
$(X_r,X_s)$ where $(X_0,X_1)$ is any random variable of
law $\Pi$. As we already said, $\Pi_{r,s}$ is an optimal
coupling of $\mu_r,\mu_s$ whose corresponding
geodesic is the restriction of $(\mu_t)$ to $[r,s]$.

Since $\Pi_{r,s}$ is optimal, according to the cyclical
monotonicity (see Lemma \ref{lemm:cyclical_monotonicity})
one has $(y_r-x_r)\cdot(y_s-x_s)\geqslant0$.

But with the above notations, one has
$y_r-x_r=u+rv$ and $y_s-x_s=u+sv$,
and we get the desired inequality. 
\end{proof}

Let us show why this Lemma implies that the only complete geodesics
are those obtained by translation. There are immediate
consequences on the rank of $\wass_2(\mR^n)$,
see Theorem \ref{enonce:flats} and Section \ref{sec:ranks}.

\begin{prop}\label{enonce:complete_geod}
Let $\mu=(\mu_t)_{t\in\mR}$ be a geodesic in $\wass_2(\mR^n)$ defined
for all times. Then there is a vector $u$
such that $\mu_t=(T_{tu})_\#\mu_0$.
\end{prop}

This result holds even if $n=1$, as stated in Proposition \ref{enonce:consequence}.

\begin{proof}
It is sufficient to find a $u$ such that 
$\mu_1=(T_u)_\#\mu_0$, since then there is
only one geodesic from $\mu_0$ to $\mu_1$.

Consider any pair of points $(x_0,x_1),(y_0,y_1)$
in the support of the coupling $\Pi$
between $\mu_0$ and $\mu_1$ that defines
the restriction of $\mu$
to $[0,1]$. Define $u=y_0-x_0$ and $v=y_1-x_1-(y_0-x_0)$.
If $v\neq0$, then there are real numbers $r<s$
such that $|u|^2+(r+s)u\cdot v+rs|v|^2<0$.
Then the coupling $\Pi_{r,s}$ between $\mu_r$ and $\mu_s$
that defines the restriction of $\mu$ to $[r,s]$,
defined as above, cannot be optimal. This is a contradiction
with the assumption that $\mu$ is a geodesic.

Therefore, for all $(x_0,x_1),(y_0,y_1)$ in the support
of $\Pi$ one has $y_0-x_0=y_1-x_1$. This amounts
to say that $\Pi$ is deterministic, given by a translation
of vector $u=y_0-x_0$.
\end{proof}

\section{Curvature}\label{sec:curvature}

Once again, this section mainly collects some facts that are 
already well-known but shall be used on the sequel.

More details on the (sectional) curvature of metric spaces are available for example
in \cite{Burago} or \cite{Jost}.
We shall consider the curvature of $\wass_2(\mR)$, in the sense of
Alexandrov. Given any three points $x,y,z$ in a geodesic metric space $X$,
there is up to congruence a unique comparison triangle $x',y',z'$ in $\mR^2$,
that is a triangle that satisfies $d(x,y)=d(x',y')$, $d(y,z)=d(y',z')$,
and $d(z,x)=d(z',x')$.

\begin{figure}[ht]\begin{center}
\input{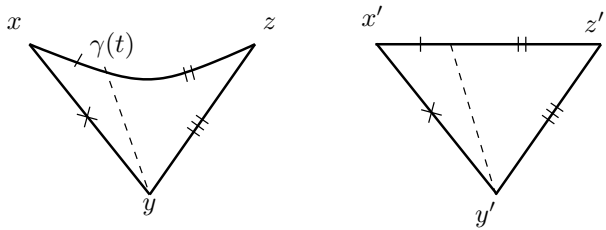}
\caption{The $\CAT(0)$ inequality: the dashed segment is shorter in the triangle $xyz$
than in the comparison triangle on the right.}\label{fig:comparison}
\end{center}\end{figure}

One says that $X$ has \emph{non-positive curvature} (in the sense of Alexandrov),
or is $\CAT(0)$,
if for all $x,y,z$ the distances between two points on sides of this triangle
is lesser than or equal to the distance between the corresponding points in the comparison
triangle, see figure \ref{fig:comparison}.

Equivalently, $X$ is $\CAT(0)$ if for any triangle $x,y,z$, any geodesic
$\gamma$ such that $\gamma(0)=x$ and $\gamma(1)=y$, and any $t\in[0,1]$,
the following inequality holds:
\begin{equation}
d^2(y,\gamma(t))\leqslant (1-t)d^2(y,\gamma(0))+td^2(y,\gamma(1))-t(1-t)t\ell(\gamma)^2
\label{eq:cat0}
\end{equation}
where $\ell(\gamma)$ denotes the length of $\gamma$, that is $d(x,z)$.

One says that $X$ has \emph{vanishing curvature}
if equality holds for all $x,y,z,\gamma,t$:
\begin{equation}
d^2(y,\gamma(t)) = (1-t)d^2(y,\gamma(0))+td^2(y,\gamma(1))-t(1-t)t\ell(\gamma)^2
\label{eq:vanishing}
\end{equation}

This is equivalent to the condition that for any triangle $x,y,z$ in $X$ and
any point $\gamma(t)$ on any geodesic segment between $x$ and $z$, the
distance between $y$ and $\gamma(t)$ is \emph{equal} to the corresponding distance
in the comparison triangle.

\begin{prop}\label{enonce:curvature}
The space $\wass_2(\mR)$ has vanishing Alexandrov curvature.
\end{prop}

\begin{proof}
It follows from the expression \eqref{eq:distance} of the distance in $\wass_2(\mR)$:
if we denote by $A,B,C$ the inverse distribution functions of the three considered
points $x,y,z\in \wass_2(\mR)$, we get:
\begin{eqnarray*}
d^2(y,\gamma(t)) &=& \int_0^1 \left(B-(1-t)A-tC\right)^2 \\
   &=& \int_0^1 \big[(1-t)^2(B-A)^2+t^2(B-C)^2 \\
   &&  + 2t(1-t)(B-A)(B-C)\big]
\end{eqnarray*}
and using that $(1-t)^2=(1-t)-t(1-t)$ and $t^2=t-t(1-t)$,
\begin{eqnarray*}
d^2(y,\gamma(t)) &=& (1-t)d^2(y,x)+td^2(y,z)
        -t(1-t)\int_0^1 \Big[(B-A)^2\\
   && +(B-C)^2-2(B-A)(B-C)\Big] \\
   &=& (1-t)d^2(y,x)+td^2(y,z)-t(1-t)d^2(x,z).
\end{eqnarray*}
\end{proof}

We shall use the vanishing curvature of $\wass_2(\mR)$ by means of the following
result, where all subsets of $X$ are assumed to be endowed with the
induced metric (that need therefore not be inner).

\begin{prop}\label{enonce:extension}
Let $X$ be a Polish uniquely geodesic space with vanishing curvature. If $Y$ is
a subset of $X$ and $C(Y)$ is the convex hull of $Y$, then any isometry
of $Y$ can be extended into an isometry of $\overline{C(Y)}$.
\end{prop}

\begin{proof}
Let $\varphi:Y\to Y$ the isometry to be extended.
Let $x,y$ be any points lying each on one geodesic segment $\gamma,\tau:[0,1]\to X$ whose
endpoints are in $Y$. Consider the unique
geodesics $\gamma',\tau'$ that satisfy $\gamma'(0)=\varphi(\gamma(0))$, $\gamma'(1)=\varphi(\gamma(1))$
$\tau'(0)=\varphi(\tau(0))$, $\tau'(1)=\varphi(\tau(1))$ and the points $x'$, $y'$ lying on them
so that $d(x',\gamma'(0))=d(x,\gamma(0))$, $d(x',\gamma'(1))=d(x,\gamma(1))$,
and the same for $y'$. This makes sense since, $\varphi$ being an isometry on $Y$,
$\gamma'$ has the length of $\gamma$ and $\tau'$ that of $\tau$. We shall prove
that $d(x',y')=d(x,y)$.

The vanishing of curvature implies that $d(x',\tau'(0))=d(x,\tau(0))$: the triangles
$\gamma(0),\gamma(1),\tau(0)$ and $\gamma'(0),\gamma'(1),\tau'(0)$ have the same
comparison triangle. Similarly $d(x',\tau'(1))=d(x,\tau(1))$. Now $x,\tau(0),\tau(1)$
and $x',\tau'(0),\tau'(1)$ have the same comparison triangle, and the vanishing curvature
assumption implies $d(x',y')=d(x,y)$.

\begin{figure}[ht]\begin{center}
\input{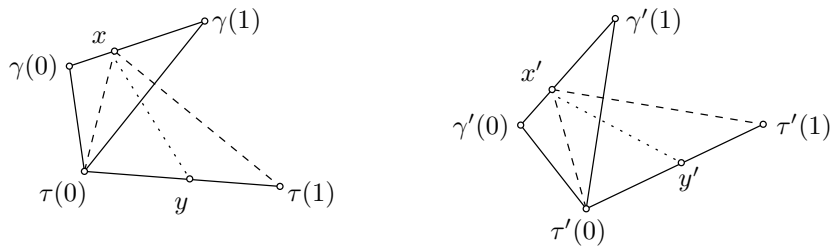}
\caption{All triangles being flat, the distance is the same between $x'$ and $y'$ and
between $x$ and $y$.}
\end{center}\end{figure}

In particular, if $x=y$ then $x'=y'$. We can thus extend $\varphi$ to the union of geodesic
segments whose endpoints are in $Y$ by mapping any such $x$ to the corresponding $x'$.
This is well-defined, and an isometry. Repeating this operation we can extend $\varphi$
into an isometry of $C(Y)$. But $X$ being complete, the continuous extension of
$\varphi$ to $\overline{C(Y)}$ is well-defined and an isometry.
\end{proof}

Note that the same result holds with the same proof when the curvature is constant but non-zero.

\subsubsection{Higher dimensional case}

Proposition \ref{enonce:curvature} does not hold in $\wass_2(\mR^n)$. In effect, there are pairs of geodesics
that meet at both endpoints (take measures whose support lie on orthogonal subspaces
of $\mR^n$). Taking a third point in one of the two geodesics,
one gets a triangle in  $\wass_2(\mR^n)$ whose comparison triangle has its three
vertices on a line. This implies that $\wass_2(\mR^n)$ is not $\CAT(0)$.
The situation is in fact worse than that: in any neighborhood $U$ of any
point of $\wass_2(\mR^n)$ one can find two different geodesics
that meet at their endpoints. One can say that this space has positive sectional
curvature at arbitrarily small scales.

\section{Isometries: the case of the line}\label{sec:isometries}

\subsection{Existence and unicity of the non-trivial isometric flow}

In this section, we prove Theorem \ref{enonce:isometries}.

Let us start with the following consequence of
Proposition \ref{enonce:consequence}.
\begin{lemm}\label{enonce:delta}
An isometry of $\wass_2(\mR)$ must globally preserve the sets $\Delta_1$
and $\Delta_2$.
\end{lemm}

\begin{figure}[t]\begin{center}
\input{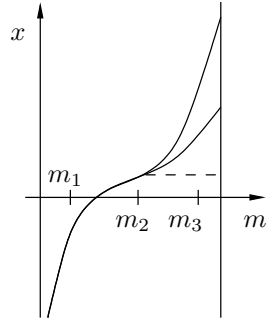}
\caption{The geodesic between these inverse distributions is defined for negative times,
more precisely until it reaches the dashed line.}\label{fig:a_geodesic}
\end{center}\end{figure}

\begin{proof}
We shall exhibit some geometric properties that characterize the points
of $\Delta_1$ and $\Delta_2$ and must be preserved by isometries.

First, according to Proposition \ref{enonce:consequence}, the points $\mu\in\Delta_1$ are the 
only ones to satisfy: every maximal geodesic ray starting at $\mu$ is complete.
Since an isometry must map a geodesic (ray, segment) to another, this
property is preserved by isometries of $\wass_2(\mR)$.

Second, let us prove that the point $\mu\in\Delta_2$ are the only ones
that satisfy: any geodesic $\mu_t$ such that $\mu=\mu_0$ and that can be extended
to a maximal interval $[T,+\infty)$ with $-\infty<T<0$, has its endpoint $\mu_T$ in 
$\Delta_1$.

This property is obviously satisfied by points of $\Delta_1$.
It is also satisfied by every points of $\Delta'_2$. Indeed, denote by $F_t$ the 
distribution function of $\mu_t$ and write $\mu=a\delta_x+(1-a)\delta_y$ where $x<y$. Then
if $\mu_1$ does not write $\mu_1=a\delta_{x_1}+(1-a)\delta_{y_1}$ with
$x_1<y_1$, either 
\begin{itemize}
\item there are two reals $b,c$ such that $a<b<c<1$ and $F_1^{-1}(b)< F_1^{-1}(c)$, 
\item there are two reals $b,c$ such that $0<b<c<a$ and $F_1^{-1}(b)< F_1^{-1}(c)$, or
\item $\mu_1$ is a Dirac mass.
\end{itemize}
In the first two cases, $\mu_t$ is not defined for $t<0$, and in the third one, it
is not defined for $t>1$.
If $\mu_1$ does write $\mu_1=a\delta_{x_1}+(1-a)\delta_{y_1}$,
then either $|y_1-x_1|=|y-x|$ and $\mu_t$ is defined for all $t$,
or $|y_1-x_1|<|y-x|$ and $\mu_t$ is only defined until a finite positive time,
or $|y_1-x_1|>|y-x|$ and $\mu_t$ is defined from a finite negative time
$T$ where $\mu_T\in\Delta_1$.

Now if $\mu\notin\Delta_2$, its inverse distribution function
$F^{-1}$ takes three different values  at some points
$m_1<m_2<m_3$. Consider the geodesic between $\mu$ and
the measure $\mu'$ whose inverse distribution function $F'^{-1}$
coincide with $F^{-1}$ on $[m_1,m_2]$ but is defined by
\[F'^{-1}(m)-F^{-1}(m_2)=2(F^{-1}(m)-F^{-1}(m_2))\]
on $[m_2,1)$ (see figure \ref{fig:a_geodesic}). Then this geodesic is defined for all positive times,
but stops at some nonpositive time $T$. Since $F^{-1}$ takes
different values at $m_2$ and $m_3$, one can extend the geodesic
for small negative times and $T<0$. But the inverse distribution function
of the endpoint $\mu_T$ must take the same values than that of $\mu$
in $m_1$ and $m_2$, thus $\mu_T\notin\Delta_1$. 
\end{proof}

Now we consider isometries of $\Delta_2$, to which all isometries of $\wass_2(\mR)$
shall be reduced.

Any point $\mu\in\Delta'_2$ writes under the form
\[\mu=\mu(x,\sigma,p)=\frac{e^{-p}}{e^{-p}+e^{p}}\delta_{x-\sigma e^p}
     +\frac{e^{p}}{e^{-p}+e^{p}}\delta_{x+\sigma e^{-p}}\]
where $x$ is its center of mass, $\sigma$ is the distance between $\mu$ and its center of mass,
and $p$ is any real number. In probabilistic terms, if $\mu$ is the law of a random variable
then $x$ is its expected value and $\sigma^2$ its variance.

\begin{lemm}\label{enonce:isom}
An isometry of $\wass_2(\mR)$ that fixes each point of $\Delta_1$ 
must restrict to $\Delta'_2$ to a map of the form:
\[\Phi(\varphi)=\mu(x,\sigma,p)\mapsto\mu(x,\sigma,\varphi(p))\]
for some $\varphi\in\isom(\mR)$. Any such map is an isometry of $\Delta_2$.
\end{lemm}

\begin{proof}
Let $\Phi$ be an isometry of $\wass_2(\mR)$ that fixes each point of $\Delta_1$.

A computation gives the following expression for the distance between two measures
in $\Delta'_2$: 
\[d^2_{\wass}(\mu(x,\sigma,p),\mu(y,\rho,q))=(x-y)^2+\sigma^2+\rho^2-2\sigma\rho e^{|p-q|}\]

Since $\Phi$ is an isometry, it preserves the center of mass and variance. The preceding
expression shows that it must preserve the euclidean distance between $p$ and $q$ 
for any two measures $\mu(x,\sigma,e),\mu(y,\rho,f)$, and that this condition is sufficient to
make $\Phi$ an isometry of $\Delta_2$.
\end{proof}

\begin{lemm}\label{enonce:semi-direct}
Let $\psi:x\to \varepsilon x+v$ and $\varphi:p\to \eta p+t$ be isometries
of $\mR$. Then 
$$\#(\psi)\Phi(\varphi)\#(\psi)^{-1}(\mu(x,\sigma,p))=\mu(x,\sigma,\eta p+\varepsilon t)$$
\end{lemm}

\begin{proof}
It follows from a direct computation:
\begin{eqnarray*}
\#(\psi)\Phi(\varphi)\#(\psi)^{-1}(\mu(x,\sigma,p))
  &=&\#(\psi)\Phi(\varphi)(\mu(\varepsilon x-\varepsilon v,\sigma,\varepsilon p))\\
  &=&\#(\psi)(\mu(\varepsilon x-\varepsilon v,\sigma,\eta\varepsilon p+t))\\
  &=&\mu(x,\sigma,\eta p+\varepsilon t)
\end{eqnarray*}
\end{proof}

We are now able to prove our main result.
\begin{proof}[Proof of Theorem \ref{enonce:isometries}]
First Lemma \ref{enonce:delta} says that any isometry
of $\wass_2(\mR)$ acts on $\Delta_1$ and $\Delta_2$.

Let $\mathscr{L}$ be $\#(\isom\mR)$ and $\mathscr{R}$ be the subset of $\isom\wass_2(\mR)$
consisting of isometries that fix $\Delta_1$ pointwise. Then $\mathscr{R}$ 
is a normal subgroup of $\isom\wass_2(\mR)$.

Let $\Psi$ be an isometry of $\wass_2(\mR)$. It acts isometrically
on $\Delta_1$, thus there is an isometry $\psi$ of $\mR$ such that
$\#(\psi)\Psi\in \mathscr{R}$. In particular, $\isom\wass_2(\mR)=\mathscr{L}\mathscr{R}$. Since
$\mathscr{L}\cap \mathscr{R}$ is reduced to the identity, we do have a semidirect product
$\isom\wass_2(\mR)=\mathscr{L} \ltimes \mathscr{R}$.

According to Proposition \ref{enonce:extension},
each map $\Phi(\varphi):\mu(x,\sigma,p)\mapsto\mu(x,\sigma,\varphi(p))$
extends into an isometry of $\overline{C(\Delta_2)}$, which
is $\wass_2(\mR)$ by Proposition \ref{enonce:convex_hull}.
We still denote by $\Phi(\varphi)$ this extension.
Proposition \ref{enonce:convex_hull} also shows that
an isometry of $\wass_2(\mR)$ is entirely determined
by its action on $\Delta_2$. The description of $\mathscr{R}$
now follows from Lemma \ref{enonce:isom}.

If $\sigma$ denotes the symmetry around $0\in\mR$, then  
$\Phi(\sigma)$ maps a measure $\mu\in \wass_2(\mR)$ to its symmetric with respect to its center
of mass, thus preserves shapes. Any other $\varphi\in\isom(\mR)$ is a translation or the composition
of $\sigma$ and a translation.

By the \emph{exotic isometry flow} of $\isom \wass_2(\mR)$ we mean the flow of isometries
$\Phi^t=\Phi(\varphi^t)$ obtained when $\varphi^t:p\to p+t$ is a translation. This flow
does not preserve shapes as is seen from its expression in $\Delta_2$.

At last, Lemma \ref{enonce:semi-direct} gives the asserted
description of the semidirect product.
\end{proof}

\subsection{Behaviour of the exotic isometry flow}

The definition of $\Phi^t$ is constructive, but not very explicit outside $\Delta_2$.
On $\Delta_2$, the flow tends to put most of the mass on the right of the center of mass,
very close to it, and send a smaller and smaller bit of mass far away on the left.

\begin{figure}[ht]\begin{center}
\input{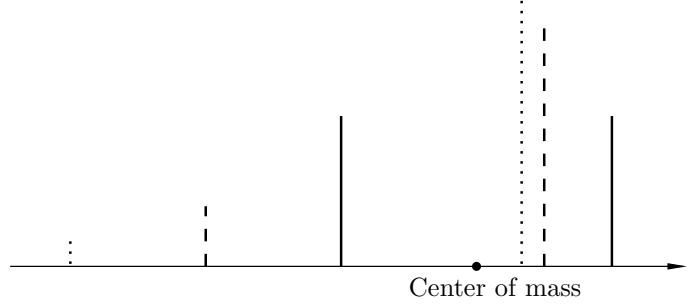}
\caption{Image of a point of $\Delta_2$ by $\Phi^2$ (dashed) and $\Phi^3$ (dotted).}
\end{center}\end{figure}

The flow $\Phi^t$ preserves $\Delta_3$ as its elements are the only ones to lie
on a geodesic segment having both endpoints in $\Delta_2$. Similarly, elements
of $\Delta_n$ are the only ones to lie on a geodesic segment having an endpoint in $\Delta_2$
and another in $\Delta_{n-1}$, therefore $\Phi^t$ preserves $\Delta_n$ for all $n$.
 
Direct computations enable one to find formulas for $\Phi^t$ on $\Delta_n$, but the expressions
one gets are not so nice. For example, if $\mu=\frac13\delta_{x_1}+\frac13\delta_{x_2}+\frac13\delta_{x_3}$
where $x_1\leqslant x_2\leqslant x_3$, then 
\begin{eqnarray*}
\Phi^t(\mu) &=& \frac{1}{1+2t^2}\ \delta_{x_1+\frac13(1-t)(x_3-x_1)+\frac13(1-t)(x_2-x_1)}\\
            && + \frac{\frac32 t^2}{\left(1+\frac12 t^2\right)\left(1+2t^2\right)}\ 
                 \delta_{x_1+\frac13(1-t)(x_3-x_1)+\frac13(1+t^{-1}+t)(x_2-x_1)}\\
            && +\frac{\frac12 t^2}{1+\frac12 t^2}\ 
                 \delta_{x_1+\frac13(1+2t^{-1})(x_3-x_1)+\frac13(1-t^{-1})(x_2-x_1)}
\end{eqnarray*}

In order to get some intuition about $\Phi^t$, let us prove the following.
\begin{prop}\label{enonce:weak_convergence}
Let $\mu$ be any point of $\wass_2(\mR)$ and $x$ its center of mass. If $t$ goes to 
$\pm\infty$, then $\Phi^t(\mu)$ converges weakly to $\delta_x$.
\end{prop}
\begin{proof}
We shall only consider the case when $t\to+\infty$ since the other one is symmetric.
Let us start with a lemma.
\begin{lemm}
If $\gamma^t$ and $\nu^t$ are in $\wass_2(\mR)$ and both converge weakly to $\delta_x$ when $t$ goes to
$+\infty$,
and if $\mu^t$ is in the geodesic segment between $\gamma^t$ and $\nu^t$ for all $t$,
then $\mu^t$ converges weakly to $\delta_x$ when $t$ goes to $+\infty$.
\end{lemm}

\begin{proof}
It is a direct consequence of the form of geodesics: if $\gamma^t$ and $\nu^t$
both charge an interval $[x-\eta,x+\eta]$ with a mass at least $1-\varepsilon$,
then $\mu^t$ must charge this interval with a mass at least $1-2\varepsilon$.
\end{proof}

Now we are able to prove the proposition on larger and larger subsets of $\wass_2(\mR)$.
First, it is obvious on $\Delta_2$. If it holds on $\Delta_n$, the preceding Lemma
together with Proposition \ref{enonce:convex} implies that it holds on $\Delta_{n+1}$.
To prove it on the whole of $\wass_2(\mR)$, the density
of the subset of $\bigcup_n\Delta_n$ consisting of measures having center of mass
$x$, and a diagonal process are sufficient.
\end{proof}

\section{Isometries: the higher-dimensional case}\label{sec:plane}

To show that the exotic isometry flow of $\wass_2(\mR)$ is exceptional,
let us consider the higher-dimensional case : there are isometries
of $\wass_2(\mR^n)$ that fix pointwise the set of Dirac masses, but there are
not so many and they preserve shapes. 

\subsection{Existence of non-trivial isometries}

We start with the existence of non-trivial isometries, that however
preserves shapes.

\begin{figure}[htp]\begin{center}
\includegraphics{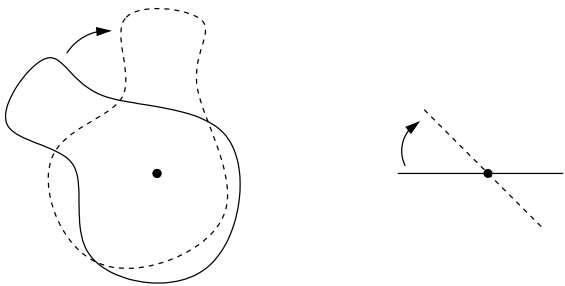}
\caption{Example of a non-trivial isometry that preserves shapes.}\label{fig:isometry_rot}
\end{center}\end{figure}

\begin{prop}
If $\varphi$ is a linear isometry of $\mR^n$, then the map
$$\Phi(\varphi):\mu\mapsto \varphi_\#(\mu-g)+g$$
where $g$ denotes both the center of mass of $\mu$ and the corresponding
translation, is an isometry of $\wass_2(\mR^n)$ (see figure \ref{fig:isometry_rot}).
\end{prop}

Note that we need $\varphi$ to be linear, thus we do not get as many non-trivial isometries as in
$\wass_2(\mR)$. Moreover, all isometries constructed this way preserve shapes.

\begin{proof}
We only need
to check the case of absolutely continuous measures $\mu,\nu$ since they
form a dense subset of $\wass_2(\mR^n)$. In that case,
there is a unique optimal coupling that is deterministic,
given by a map $T:\mR^n\to\mR^n$ such that $T_\#(\mu)=\nu$.
Denote by $g$ and $h$ the centers of mass of $\mu$ and $\nu$.
Let us show that there is a good coupling
between $\mu'=\Phi(\varphi)(\mu)$ and $\nu'=\Phi(\varphi)(\nu)$.

Let $T'$ be the map defined by
$$T'(\varphi(x-g)+g)=\varphi(Tx-h)+h$$
By construction $T'_\#(\mu')=\nu'$. Moreover the cost of the coupling
$(\id\times T')_\#\mu'$ is
\begin{eqnarray*}
A&=&\int_{\mR^n} |\varphi(x-g)+g-\varphi(Tx-h)-h|^2 \mu(dx)\\
  &=& \int_{\mR^n} \Big(|Tx-x|^2+2|g-h|^2+2(x-Tx)\cdot(h-g)+2\varphi(x-Tx)\cdot(g-h)\\
  && \qquad+2\varphi(h-g)\cdot(g-h)\Big) \mu(dx)\\
  &=& d_{\wass}^2(\mu,\nu)+2|g-h|^2+2(g-h)\cdot(h-g)+2\varphi(g-h)\cdot(g-h)\\
  && \qquad +2\varphi(h-g)\cdot(g-h) \\
  &=& d_{\wass}^2(\mu,\nu)
\end{eqnarray*}
This shows that $\mu'$ and $\nu'$ are at distance at most
$d_{\wass}(\mu,\nu)$. Applying the same reasoning to $\Phi(\varphi)^{-1}$,
we get that $d_{\wass}(\mu',\nu')=d_{\wass}(\mu,\nu)$ and $\Phi(\varphi)$
is an isometry.
\end{proof}

\subsection{Semidirect product decomposition}

The Dirac masses are the only measures
such that any geodesic issued from them can be extended for all
times (given any other measure, the geodesic pointing to any Dirac mass
cannot be extended past it, see in Section \ref{sec:Rn_bases}
the paragraph on dilations). As a consequence,
an isometry of $\wass_2(\mR^n)$ must globally preserve
the set of Dirac masses

As in the case of the line, if we let $\mathscr{L}=\#\isom(\mR^n)$
and $\mathscr{R}$ be the set of isometries of $\wass_2(\mR^n)$
that fix each Dirac mass, then
$\isom\wass_2(\mR^n)=\mathscr{L}\ltimes \mathscr{R}$.
We proved above that $\mathscr{R}$ contains a copy of $\mathrm{O}(n)$,
and is in particular non-trivial.

Moreover, if one conjugates a $\Phi(\varphi)$ by some $\#(\psi)$,
where $\varphi\in \mathrm{O}(n)$ and $\psi\in\isom\mR^n$,
one gets the map $\Phi(\vec{\psi}\varphi\vec{\psi}^{-1})$ (it is
sufficient to check this on some measure $\mu$ in 
the easy cases when $\psi$ is a translation or fixes the center
of mass of $\mu$).
Therefore, the action of $\mathscr{L}$ on $\mathrm{O}(n)\subset\mathscr{R}$
in the semidirect product is as asserted in Theorem \ref{enonce:isom_highdim}.

To deduce Theorem \ref{enonce:isom_highdim}, 
we are thus left with proving that an isometry that fixes pointwise all
Dirac masses must be of the form $\Phi(\varphi)$ for some
$\varphi\in \mathrm{O}(n)$. 

\subsection{Measures supported on subspaces}

The following lemma will be used to prove that isometries of
$\wass_2(\mR^n)$ must preserve the property of being supported on a
proper subspace.

\begin{lemm}\label{lemm:orthogonal_measures}
Let $\mu,\nu\in\wass_2(\mR^n)$, denote by $g,h$ their centers of
mass and let $\sigma=d_{\wass}(\mu,\delta_g)$ and $\rho=d_{\wass}(\nu,\delta_h)$.
The equality 
\begin{equation}
d_{\wass}^2(\mu,\nu)=d^2(g,h)+\sigma^2+\rho^2
\label{eq:indep}
\end{equation}
holds if and only if there are two orthogonal affine subspaces $L$ and $M$ such that
$\mu\in \wass_2(L)$ and $\nu\in \wass_2(M)$.
\end{lemm}

\begin{proof}
Let us first prove that $d^2(g,h)+\sigma^2+\rho^2$ is the cost $B$ of the independent
coupling $\Pi=\mu\otimes\nu$:
\begin{eqnarray*}
B &:=& \int_{\mR^n\times\mR^n}d^2(x,y)\Pi(dxdy) \\
  &=& \int_{\mR^n\times\mR^n} |(x-g)-(y-h)+(g-h)|^2 \Pi(dxdy)\\
  &=&\sigma^2+\rho^2+d^2(g,h)-2\Big(\int_{\mR^n}(x-g)\mu(dx)\Big)\cdot\Big(\int_{\mR^n}(y-h)\nu(dy)\Big)\\
  &&\qquad +2(g-h)\cdot\Big(\int_{\mR^n}(x-g)\mu(dx)\Big) -2(g-h)\cdot\Big(\int_{\mR^n}(y-h)\nu(dy)\Big)\\
  &=& \sigma^2+\rho^2+d^2(g,h)
\end{eqnarray*}
since by definition $g=\int x\mu(dx)$ and $h=\int y\mu(dy)$.

As a consequence, \eqref{eq:indep} holds if and only if the independent coupling
is optimal. 

If $\mu$ has two point $x,y$ in its support and $\nu$ has two points $z,t$ in its support
such that $(xy)$ is not orthogonal to $(zt)$, then either $(x-y)\cdot(z-t)<0$ or
$(x-y)\cdot(t-z)<0$. Then by cyclical monotonicity (see Lemma \ref{lemm:cyclical_monotonicity})
the support of an optimal coupling cannot contain $(x,z)$ and $(y,t)$ (in the first case)
or $(z,x)$ and $(y,t)$ (in the second case) and thus cannot be the independent coupling.

As a consequence, \eqref{eq:indep} holds if and only if $\mu$ and $\nu$ are supported on two orthogonal
affine subspaces.
\end{proof}

\begin{lemm}\label{lemm:hyperplane_supported}
Isometries of $\wass_2(\mR^n)$ send hyperplane supported-measures on
hyperplane supported measures. Moreover, if two measures are supported
on parallel hyperplanes, then their images by any isometry are supported by
parallel hyperplanes
\end{lemm}

\begin{proof}
Let $\mu\in\wass_2(\mR^n)$ be supported by some hyperplane $H$ and $\Phi$
be an isometry of $\wass_2(\mR^n)$. Let $\nu$ be any measure that is supported on
a line orthogonal to $H$, and that is not a Dirac mass. Then \eqref{eq:indep} holds
(whith the same notation as above). 

Let $\mu'$ and $\nu'$ denote the images of $\mu$ and $\nu$ by $\Phi$.
We know that $\Phi$ must map $\delta_g$ and $\delta_h$
to Dirac masses $\delta_{g'}$ and $\delta_{h'}$. Since the center of mass of an element of $\wass_2(\mR^n)$ is
uniquely defined as its projection on the set of Dirac masses, $g'$ is the center of mass
of $\mu'$ and $h'$ is that of $\nu'$. We get, denoting by $\sigma'$ and $\rho'$
the distances of $\mu'$ and $\nu'$ to their centers of mass,
\begin{eqnarray*}
d_{\wass}^2(\mu',\nu') &=& d_{\wass}^2(\mu,\nu)\\
                       &=& d^2(g,h)+\sigma^2+\rho^2 \\
                       &=& d^2(g',h')+\sigma'^2+\rho'^2
\end{eqnarray*}
which implies that $\mu'$ and $\nu'$ are supported on orthogonal subspaces $L$ and $M$ of $\mR^n$.
Since $\nu'$ is not a Dirac mass, $M$ is not a point and $L$ is contained in some hyperplane $H'$.

Moreover, if $\mu_1$ is another measure supported on a hyperplane parallel to $H$, then
its image $\mu_1'$ is supported on some subspace orthogonal to $M$. It follows that we can find
a hyperplane parallel to $H'$ that contains the support of $\mu_1'$.
\end{proof}

We are know equipped to end the proof of Theorem \ref{enonce:isom_highdim} by induction
on the dimension.

\subsection{Non-existence of exotic isometries: case of the plane}

Given a line $L\subset\mR^2$,
denote by $\wass_2(L)$ the subset of $\wass_2(\mR^2)$
consisting of all measures whose support is a subset of
$L$. An optimal coupling between two points of $\wass_2(L)$
must have its support in $L\times L$, thus $\wass_2(L)$
endowed with the restriction of the distance of $\wass_2(\mR^2)$
is isometric to $\wass_2(\mR)$. More precisely, given any
isometry $\psi : \mR\to L$, we get an isometry 
$\psi_\#:\wass_2(\mR)\to \wass_2(L)$.

Lemma \ref{lemm:hyperplane_supported} ensures
that any isometry $\Phi$
maps a line-supported measure to a line-supported measure, and that
the various measures in  $\wass_2(L)$ are mapped to measures
supported on parallel lines. 

We assume from now on that $\Phi$ fixes each
Dirac mass and, up to composing it with some
$\Phi(\varphi)$, that it preserves globally
$\wass_2(L)$ for some $L$ (the axis $\mR\times\{0\}$ say). We can moreover
assume that its restriction to $\wass_2(L)$ is a $\Phi^t$ for some $t$.

\begin{lemm}
Let $\Phi$ be an isometry of $\wass_2(\mR^2)$ that fixes Dirac masses,
preserves globally $\wass_2(L)$ and such that its restriction to this subspace
is the time $t$ of the exotic isometric flow.

Then $t$ must be $0$ and up to composing with some
$\Phi(\varphi)$, we can assume that $\Phi$ preserves $\wass_2(M)$ for all line $M$.
\end{lemm}

\begin{proof}
We identify a measure $\mu\in\wass_2(\mR)$ with its image
by the usual embedding that identifies $\mR$ with the axis $L$.
We denote by $\theta L$ the rotate of
$L$ by an angle $\theta$ around the origin.

Denote by $\mu(x,\sigma,p,\theta)$ the combination of two Dirac masses
that is the image of $\mu(x,\sigma,p)$ if $\theta=0$, and its rotate
around $x$ by an angle $\theta$ otherwise. If $\theta\leqslant\pi/2$, one gets
\begin{equation}
d_{\wass}^2(\mu(0,1,p,0),\mu(0,1,q,\theta))=2-2e^{|p-q|}\cos\theta
\label{eq:theta}
\end{equation}

This shows in particular that the measures supported on $\theta L$
and with center of mass $0$ must be mapped to measures supported
on $\pm \theta L$. Up to composing with $\Phi(\varphi)$ where
$\varphi$ is the orthogonal symmetry with respect to $L$, we can assume
that the measures supported on $\frac{\pi}3 L$ and with center of mass $0$
are mapped to measures supported on $\frac{\pi}3 L$.

Then the measures supported on $\theta L$
and with center of mass $0$ must be mapped to measures supported
on a line that crosses $L$ and $\frac{\pi}3 L$ with angles
$\pm \theta$ and $\pm(\theta-\frac{\pi}3)$. They are therefore mapped
to measures supported on $\theta L$.

Using the same argument and Lemma \ref{lemm:hyperplane_supported}, 
we get that $\Phi$ must preserve $\wass_2(M)$ for all lines $M$.

Moreover, from equation \eqref{eq:theta} we deduce that
if the restriction of $\Phi$ to $\wass_2(L)$
is the time $t$ of the exotic isometric flow, then for all $\theta<\pi/2$
its restriction to $\theta L$ also is. But applying the same reasoning
to $\frac{\pi}3 L$ and $\frac{2\pi}3 L$, then to $\frac{2\pi}3 L$
and $\pi L=L$ we see that the restriction of $\Phi$ to
$\wass_2(L)$ \emph{with the reversed orientation} must be the time
$t$ of the exotic isometric flow. This implies $t=-t$ thus $t=0$.
\end{proof}

The case $n=2$ of Theorem \ref{enonce:isom_highdim} is
now reduced to the following.
\begin{lemm}
If an isometry $\Phi$ of $\wass_2(\mR^2)$ fixes pointwise the set of line-supported
measures, then it must be the identity.
\end{lemm}

\begin{proof}
This is a consequence of Radon's theorem, which asserts that
a function (compactly supported and smooth, say) in $\mR^n$ is characterized by its
integrals along all hyperplanes \cite{Radon} (also see \cite{Helgason}).

Given $\mu\in\wass_2(\mR^2)$,  one can determine
by purely metric means its orthogonal projection on any fixed line $L$: it is
its metric projection, that is the unique $\nu$ supported on $L$
that minimizes the distance $d_{\wass}(\mu,\nu)$.

Now if $\mu$ has smooth density and is compactly supported,
then the integral of its density along any line $L$ is
exactly the density at point $M\cap L$ of its projection onto any line $M$
orthogonal to $L$.

Therefore, $\Phi$ must fix every measure $\mu\in\wass_2(\mR^2)$
that has a smooth density and is compactly supported.
They form a dense set of $\wass_2(\mR^2)$ thus $\Phi$ must be the identity.
\end{proof}

\subsection{Non-existence of exotic isometries: general case}

We end the proof of Theorem \ref{enonce:isom_highdim} by an induction on the dimension.
There is nothing new compared to the case of the plane, so we stay sketchy.

Let $\Phi$ be an isometry of $\wass_2(\mR^n)$ that
fixes pointwise the Dirac masses. It
must map every hyperplane-supported measure to
a hyperplane-supported measure. Using a non-trivial isometry,
we can assume that for some hyperplane $L$,
$\Phi$ globally preserves the set $\wass_2(L)$ of measures
supported on $L$. 
Thanks to the induction hypothesis, we can compose
$\Phi$ with another non-trivial isometry to ensure that
$\Phi$ fixes $\wass_2(L)$ pointwise.

Let $\mu$ be a measure supported on some hyperplane $M\neq L$.
Let $M'$ be a hyperplane supporting $\Phi(\mu)$. Then as in the
case of the plane, it is easy to show that the dihedral
angle of $(L,M)$ equals that of $(L,M')$. Moreover,
all measures supported on $L\cap M$ are fixed by $\Phi$,
and we conclude that $M'=M$ (up to composition with $\Phi(\varphi)$
where $\varphi$ is the orthogonal symmetry with respect to $L$).

The same argument shows that $\Phi$ preserves $\wass_2(M)$ for all hyperplanes
$M$.

A measure supported on $M$
is determined, if its dihedral angle with $L$ different from $\pi/2$,
by its orthogonal projection onto $L$. since $\Phi$ fixes $\wass_2(L)$
pointwise, it must fix $\wass_2(M)$ pointwise as well.
When $M$ is orthogonal to $L$, the use of a third hyperplane not orthogonal
to $M$ nor $L$ yields the same conclusion.

Now that we know that $\Phi$ fixes every hyperplane-supported
measure, we can use the Radon Theorem to conclude that it is the identity.

\section{Ranks}\label{sec:ranks}

One usually defines the \emph{rank} of a metric space $X$ as the supremum
of the set of positive integers $k$ such that there is an isometric embedding of $\mR^k$
into $X$.

As a consequence of Proposition \ref{enonce:complete_geod}, we get the following result
announced in the introduction.
\begin{theo}
The space $\wass_2(\mR^n)$ has rank $n$.
\end{theo}

\begin{proof}
An isometric embedding $e:\mR^{n+1}\to \wass_2(\mR^n)$ must map a geodesic to a geodesic, since they are
precisely those curves $\gamma$ satisfying 
\[d(\gamma(t),\gamma(s))=v|t-s|\]
 for some constant
$v$. The union of complete geodesics through any point $\mu$ in the image of $e$
 would contain a copy of $\mR^{n+1}$, but Proposition \ref{enonce:complete_geod}
shows that this union is isometric to $\mR^n$.
\end{proof}

However, one can define less restrictive notions of rank as follows.

\begin{defi}
Let $X$ be a Polish space. 
The \emph{semi-global rank} of $X$
is defined as the supremum of the set of positive integers $k$ such that
for all $r\in\mR^+$, there is an isometric embedding
of the ball of radius $r$ of $\mR^k$ into $X$.

The \emph{loose rank} of $X$ is defined as the supremum of the set of positive
integers $k$ such that there is a quasi-isometric embedding of $\mZ^k$ into $X$.
\end{defi}

Let us recall that a map $f:Y\to X$ is said to be a quasi-isometric embedding if there
are constants $C>1,D>0$ such that for all $y,z\in Y$ the following holds :
\[C^{-1}d(y,z)-D\leqslant d(f(y),f(z)) \leqslant Cd(y,z)+D.\]

The notion of loose rank is relevant in a large class of metric spaces, including
discrete spaces (the Gordian space \cite{Gambaudo-Ghys}, or the Cayley
graph of a finitely presented group for example). We chose not to call it ``coarse
rank'' due to the previous use of this term by Kapovich, Kleiner and Leeb.

The semi-global rank is motivated by the following simple result.
\begin{prop}\label{enonce:rank_hyperbolicity}
A geodesic space $X$ that has semi-global rank at least $2$ is not 
$\delta$-hyperbolic.
\end{prop}

\begin{proof}
Since $X$ contains euclidean disks of arbitrary radius, it also contains
euclidean equilateral triangles of arbitrary diameter. In such a triangle, the
maximal distance between a point of an edge and the other edges is proportional
to the diameter, thus is unbounded in $X$.
\end{proof}

\begin{prop}\label{enonce:ranks_line}
The semi-global rank and the loose rank of $\wass_2(\mR)$ are infinite.
\end{prop}

\begin{proof}
Consider the subset $\mR^k_\leqslant=\{(x_1,\ldots,x_k) ; x_1\leqslant x_2\leqslant\cdots\leqslant x_k\}$
of $\mR^k$. It is a closed, convex cone.

Moreover the map
\begin{eqnarray*}
\mR^k_\leqslant &\to& \wass_2(\mR) \\
(x_1,\ldots,x_k)&\mapsto&\sum\frac{1}{k}\delta_{x_i}
\end{eqnarray*}
is an isometric embedding.

Since $\mR^k_\leqslant$ contains arbitrarily large balls, $\wass_2(\mR)$ must have
infinite semi-global rank.

Moreover, since $\mR^k_\leqslant$ is a convex cone of non-empty interior, it
contains a circular cone. Such a circular cone is conjugate by a linear (and thus bi-Lipschitz)
map to the cone 
\[{\mathscr C}=\{x_1^2=\sum_{i\geqslant 2} x_i^2\}.\]
Now the vertical projection from $\{x_1=0\}$ to ${\mathscr C}$ is bi-Lipschitz. There is 
therefore a bi-Lipschitz embedding of $\mR^{k-1}$ in $\wass_2(\mR)$ and, \textit{a fortiori},
a quasi-isometric embedding of $\mZ^{k-1}$. Therefore $\wass_2(\mR)$ has infinite loose
rank.
\end{proof}

\subsection{Ranks of other spaces}

The ranks of $\wass_2(\mR)$ have an influence on those of many spaces due to the following
lemma.

\begin{lemm}
If $X$ and $Y$ are Polish geodesic spaces, any isometric embedding
$\varphi:X\to Y$ induces an isometric embedding 
$\varphi_\#:\wass_2(X)\to \wass_2(Y)$.
\end{lemm}
As usual, $\varphi_\#$ is defined by: $\varphi_\#\mu(A)=\mu(\varphi^{-1}(A))$ for all measurable
$A\subset Y$.

\begin{proof}
Since $\varphi$ is isometric, for any $\mu\in \wass_2(X)$, $\varphi_\#\mu$ is in $\wass_2(Y)$.
Moreover any optimal
transportation plan in $X$ is mapped to an optimal transportation plan in $Y$ (note that a coupling
between two measures with support in $\varphi(X)$ must have its support contained
in $\varphi(X)\times\varphi(X)$). Integrating the equality $d(\varphi(x),\varphi(y))=d(x,y)$
yields the desired result.
\end{proof}

\begin{coro}
If $X$ is a Polish geodesic space that contains a complete geodesic, then
$\wass_2(X)$ has infinite semi-global rank and infinite loose rank. As a consequence,
$\wass_2(X)$ is not $\delta$-hyperbolic.
\end{coro}

\begin{proof}
Follows from the preceding Lemma, Proposition \ref{enonce:ranks_line} and Proposition 
\ref{enonce:rank_hyperbolicity}.
\end{proof}

This obviously applies to $\wass_2(\mR^n)$.

One could hope that in Hadamard spaces, the projection to the center of mass
\[P:\wass_2(X)\to X\]
 could give
a higher bound on the rank of $\wass_2(X)$ by means of that of $X$. However, $P$ need not
map a geodesic on a geodesic. For example, if one consider on the real hyperbolic plane
${\mathbb R}{\mathrm H}^2$ the measures $\mu_t=1/2\delta_p+1/2\delta_{\gamma(t)}$
where $p$ is a fixed point and $\gamma(t)$ is a geodesic, then $\mu_t$ is a geodesic of
$\wass_2({\mathbb R}{\mathrm H}^2)$ that is mapped by $P$ to a curve with the same endpoints
than $\gamma$, but is different from it. Therefore, it cannot be a geodesic.

\begin{figure}[ht]\begin{center}
\input{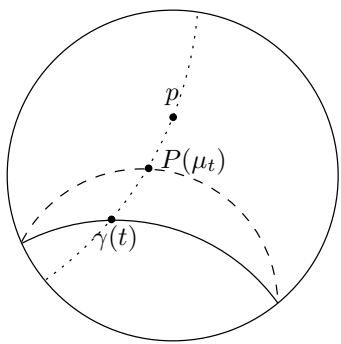}
\caption{The projection $P$ maps a geodesic of $\wass_2({\mathbb R}{\mathrm H}^2)$
to a non-geodesic curve (dashed) in ${\mathbb R}{\mathrm H}^2$.}
\end{center}\end{figure}

\section{Open problems}\label{sec:open}

Since the higher-dimensional Euclidean spaces are more rigid (has few non-trivial isometries) than the
line, we expect other spaces to be even more rigid. 
\begin{ques}
Does it exists a Polish (or Hadamard) space $X\neq\mR$ such that $\wass_2(X)$ admits
exotic isometries?

Does it exists a Polish (or Hadamard) space $X\neq\mR^n$ such that $\wass_2(X)$ admits
non-trivial isometries?
\end{ques}
In any Hadamard space $X$, isometries
of $\wass_2(X)$ must preserve the set of Dirac masses (the proof is the same
than in $\mR$), and this fact could help get a grip on the problem in this case.

For general spaces, even the following seems not obvious.
\begin{ques}
Does it exists a Polish space $X$ whose Wasserstein space $\wass_2(X)$
possess an isometry that does not preserve the set of Dirac masses ?
\end{ques}

Last, when $X$ is Hadamard, one could
hope to use the projection $P$ to link the rank of $\wass_2(X)$ to
the loose rank of $X$.

\begin{ques}
If $X$ is a Hadamard space, is the loose rank of $X$ an upper bound
for the rank of $\wass_2(X)$ ?
\end{ques}

\bibliographystyle{smfplain}
\bibliography{biblio.bib}

\providecommand{\bysame}{\leavevmode ---\ }
\providecommand{\og}{``}
\providecommand{\fg}{''}
\providecommand{\smfandname}{et}
\providecommand{\smfedsname}{\'eds.}
\providecommand{\smfedname}{\'ed.}
\providecommand{\smfmastersthesisname}{M\'emoire}
\providecommand{\smfphdthesisname}{Th\`ese}
\begin{thebibliography}{10}

\bibitem{Ambrosio-Gigli}
{\scshape L.~Ambrosio {\normalfont \smfandname} N.~Gigli} -- {\og Construction
  of the parallel transport in the wasserstein space\fg}, \emph{Methods Appl.
  Anal.} \textbf{15} (2008), no.~1, p.~1--30.

\bibitem{Ambrosio-Gigli-Savare}
{\scshape L.~Ambrosio, N.~Gigli {\normalfont \smfandname} G.~Savar{\'e}} --
  \emph{Gradient flows in metric spaces and in the space of probability
  measures}, Lectures in Mathematics ETH Z\"urich, Birkh\"auser Verlag, Basel,
  2005.

\bibitem{Brenier1}
{\scshape Y.~Brenier} -- {\og D\'ecomposition polaire et r\'earrangement
  monotone des champs de vecteurs\fg}, \emph{C. R. Acad. Sci. Paris S\'er. I
  Math.} \textbf{305} (1987), no.~19, p.~805--808.

\bibitem{Brenier2}
\bysame , {\og Polar factorization and monotone rearrangement of vector-valued
  functions\fg}, \emph{Comm. Pure Appl. Math.} \textbf{44} (1991), no.~4,
  p.~375--417.

\bibitem{Bridson-Haeffliger}
{\scshape M.~R. Bridson {\normalfont \smfandname} A.~Haefliger} -- \emph{Metric
  spaces of non-positive curvature}, Grundlehren der Mathematischen
  Wissenschaften [Fundamental Principles of Mathematical Sciences], vol. 319,
  Springer-Verlag, Berlin, 1999.

\bibitem{Burago}
{\scshape D.~Burago, Y.~Burago {\normalfont \smfandname} S.~Ivanov} -- \emph{A
  course in metric geometry}, Graduate Studies in Mathematics, vol.~33,
  American Mathematical Society, Providence, RI, 2001.

\bibitem{Gambaudo-Ghys}
{\scshape J.-M. Gambaudo {\normalfont \smfandname} {\'E}.~Ghys} -- {\og Braids
  and signatures\fg}, \emph{Bull. Soc. Math. France} \textbf{133} (2005),
  no.~4, p.~541--579.

\bibitem{Helgason}
{\scshape S.~Helgason} -- \emph{The {R}adon transform}, Progress in
  Mathematics, vol.~5, Birkh\"auser Boston, Mass., 1980.

\bibitem{Jost}
{\scshape J.~Jost} -- \emph{Nonpositive curvature: geometric and analytic
  aspects}, Lectures in Mathematics ETH Z\"urich, Birkh\"auser Verlag, Basel,
  1997.

\bibitem{Juillet}
{\scshape N.~Juillet} -- {\og Optimal transport and geometric analysis in
  heisenberg groups\fg}, \smfphdthesisname, Universit{\'e} Grenoble 1 et
  Rheinische Friedrich-Wilhelms-Universit{\"a}t Bonn, December 2008.

\bibitem{KS1}
{\scshape M.~Knott {\normalfont \smfandname} C.~S. Smith} -- {\og On the
  optimal mapping of distributions\fg}, \emph{J. Optim. Theory Appl.}
  \textbf{43} (1984), no.~1, p.~39--49.

\bibitem{Villani-Lott}
{\scshape J.~Lott {\normalfont \smfandname} C.~Villani} -- {\og Ricci curvature
  for metric-measure spaces via optimal transport\fg}, \emph{Ann. of Math. (2)}
  \textbf{169} (2009), no.~3, p.~903--991.

\bibitem{Lott}
{\scshape J.~Lott} -- {\og Some geometric calculations on {W}asserstein
  space\fg}, \emph{Comm. Math. Phys.} \textbf{277} (2008), no.~2, p.~423--437.

\bibitem{Radon}
{\scshape J.~Radon} -- {\og {\"U}ber die bestimmung von funktionen durch ihre
  integralwerte l{\"a}ngs gewisser mannigfaltigkeiten\fg}, \emph{Ber. Verh.
  S{\"a}chs. Akad. Wiss. Leipzig, Math-Nat. Kl.} \textbf{69} (1917),
  p.~262--277.

\bibitem{vonRenesse-Sturm}
{\scshape M.-K. von Renesse {\normalfont \smfandname} K.-T. Sturm} -- {\og
  Transport inequalities, gradient estimates, entropy, and {R}icci
  curvature\fg}, \emph{Comm. Pure Appl. Math.} \textbf{58} (2005), no.~7,
  p.~923--940.

\bibitem{KS2}
{\scshape C.~S. Smith {\normalfont \smfandname} M.~Knott} -- {\og Note on the
  optimal transportation of distributions\fg}, \emph{J. Optim. Theory Appl.}
  \textbf{52} (1987), no.~2, p.~323--329.

\bibitem{Sturm}
{\scshape K.-T. Sturm} -- {\og On the geometry of metric measure spaces. {I},
  {II}\fg}, \emph{Acta Math.} \textbf{196} (2006), no.~1, p.~65--131, 133--177.

\bibitem{Takatsu}
{\scshape A.~Takatsu} -- {\og On wasserstein geometry of the space of gaussian
  measures\fg}, arXiv:0801.2250, 2008.

\bibitem{Villani}
{\scshape C.~Villani} -- \emph{Topics in optimal transportation}, Graduate
  Studies in Mathematics, vol.~58, American Mathematical Society, Providence,
  RI, 2003.

\bibitem{Villani2}
\bysame , \emph{Optimal transport old and new}, Grundlehren der Mathematischen
  Wissenschaften [Fundamental Principles of Mathematical Sciences], vol. 338,
  Springer-Verlag, Berlin, 2009.

\end{thebibliography}

\end{document}